\documentclass{amsart}
\usepackage{graphicx} 
\usepackage{amssymb}
\newtheorem{theorem}{Theorem}
\newtheorem{proposition}{Proposition}
\newtheorem{lemma}{Lemma}

\newtheorem{corollary}{Corollary}

\newcommand{\N} {{\rm I}\!{\rm N}}
\newcommand{\R} {{\rm I}\!{\rm R}}
\newcommand{\Z} {{\rm Z}\!\!{\rm Z}}

\input tcilatex

\def\endproof{{\hfill{\vrule height 1.5ex width 1.5ex depth -.1ex}}}
\def\upp#1{\leavevmode \raise.50ex\hbox{#1}}
\def\up#1{\leavevmode \raise.38ex\hbox{#1}}

\def\grau#1{{\hbox{{\hbox{#1}}\kern-.25em\accent'27$\!$}}}

\def\R#1{{\hbox{\bf {R}$^{#1}$}}}

\def\lt{\longrightarrow}
\def\leva{\longmapsto}

\def\h0{{\hbox{O}}}

\countdef\balance=10
\countdef\vale=1

\def\lee{ < \>
{\hbox{\kern-.45em\llap{\lower1.2ex\hbox{$ = $}}}}
{\hbox{\kern.45em{}}}}

\def\gee{ > \>
{\hbox{\kern-.45em\llap{\lower1.2ex\hbox{$ = $}}}}
{\hbox{\kern.45em{}}}}

\def\trv{ \ \big\vert
{\hbox{\kern.158em\llap{\lower-0.2ex\hbox{$ \cap $}}}}
{\hbox{\kern-.158em{}}}\ }

\def\subsetneq{\subseteq
{\hbox{\kern-.45em\llap{\lower.5ex\hbox{$\scriptstyle{/}$}}}}
{\hbox{\kern.45em{}}}}

\def\abs#1{\vert {#1} \vert}

\def\nrm#1{\Vert {#1} \Vert}

\def\fnc#1#2#3{ \ifnum\balance=0  {#1}\colon{#2} \leva {#3} 
                \else
                \ifnum\balance=1 {#1}\colon{#2} \lt {#3} 
                \else {#1}\colon {#2} & \lt{#3} \fi\fi }

\def\vai#1#2{ \ifnum\balance=0  {#1} & \leva {#2} \else
               {#1} \leva  {#2} \fi }

\begin{document}

\title{ Transitivity and rotation sets with nonempty interior for homeomorphisms of the 2-Torus }
\author{F\'abio Armando Tal}
\address{Instituto de Matem\'atica e Estat\'\i stica, Universidade de S\~ao Paulo,  Rua do Mat\~ao 1010, Cidade Universit\'aria, 05508-090 S\~ao Paulo, SP, Brazil}
\email{fabiotal@ime.usp.br}
\thanks{Supported by CNPq grant 304360/05-8}

\begin{abstract}
We show that, if $f$ is a homeomorphism of the 2--torus isotopic to the identity, and its lift $\widetilde f$ is transitive, or even if it is transitive outside of the lift of the elliptic islands, then $(0,0)$ is in the interior of the rotation set of $\widetilde f.$ This proves a particular case of Boyland´s conjecture.

\end{abstract}

\maketitle

\section{Introduction}


In this work we study homeomorphisms of the 2 torus isotopic to the identity. Let $f$ be such a homeomorphism, and let $\widetilde{f}$ be a fixed lift of $f$ to the plane, a homeomorphism of $\R{2}$ that satisfies
$$\widetilde{f}(x+1,y+1)=\widetilde{f}(x,y)+(1,1),$$
and $\pi(\widetilde{f}(x))=f(\pi(x)),$ where $\pi:\R{2}\to T^2$ is the covering map.

Given $f$ and $\widetilde{f}$ we define, following \cite{MZ}and \cite{franks}, the rotation set of $\widetilde{f}$ as the set of accumulation points of the subset of $\R{2}$
$$\left\{\frac{\widetilde{f}^n(x)-x}{n}\mid x\in \R{2},\,n\in\Z^{+} \right\},$$
which is a compact convex set (see \cite{MZ}).

One important question is to characterize for which homeomorphisms of the torus the rotation set of its lifts has non-empty interior. There are several consequences associated with this property. In \cite{LM}, it is shown that if a lift of $f$ has a rotation set with non-empty interior, then $f$ has positive topological entropy. In \cite{MZ} it is shown that if a point $(\frac{p_1}{q}, {\frac{p_2}{q}})\in Q^2$ is in the interior of $\rho(\widetilde{f}),$  then there exists a point $x$ in the plane such that $\widetilde{f}^q(x)=x+(p_1,p_2),$ while the results from \cite{Jaeger} and \cite{MZ2} show that, if $\alpha$ is a point in the interior of  $\rho(\widetilde{f}),$ with both coordinates irrational, then there exists a minimal set $K_{\alpha}$ such that the restriction of $f$ to $K_{\alpha}$ is regularly semi-conjugate to the translation by $\alpha.$  

While the dynamics tied to the interior of the rotation set is somewhat well understood, a description of the dynamics implied by the extremal points of the rotation set when it has interior is somewhat lacking, even for the cases where $f$ is area preserving. One standing conjecture in this direction is the Boyland conjecture, explained below:

Given an area preserving homeomorphism of the torus $f$ and a fixed lift $\widetilde f,$ one can consider the rotation vector of the Lebesgue measure $\lambda$ for $\widetilde{f},$
$$\rho_{\lambda}(\widetilde{f}) =\int_{[0,1]\times[0,1]}\widetilde{f}(x)-x d\lambda.$$ An homeomorphisms $f$ isotopic to the identity such that $\rho_{\lambda}(\widetilde{f})=(0,0)$ is usually called in the literature an irrotational homeomorphism. Whenever we speak of an irrotational homeomorphisms, an specific lift for $f$ is fixed. 

Boyland's conjecture for the torus claims that, if $f$ is an irrotational homeomorphisms and the interior of $\rho(\widetilde{f})$ is non-empty, then $(0,0)$ is in the interior of the rotation set.

Another important dynamical property of homeomorphisms which has received attention recently is transitivity. In \cite{BC} it was shown that $\mathcal{C}^1$ generically, an area preserving diffeomorphism of the torus is transitive. In \cite{AP}, the authors showed that, $\mathcal{C}^0$ generically, if $f$ is irrotational, then $\widetilde{f}$ is transitive in the whole plane. The work \cite{BC} suggests this is also true $\mathcal{C}^1$ generically.

In \cite{AT1} it is shown that, if $g$ is a homeomorphism of the closed annulus whose lift $\widetilde{g}$ to the  strip $\R{}\times[0,1]$ is transitive, and if there are no fixed point in the boundary of the annulus, then $0$ is in the interior of the rotation interval of $\widetilde{g}.$ Some of the ideas presented in this note follow from that paper. 

In \cite{nancy} the authors study the relationship between transitivity, entropy and rotation set, showing that if $f$ is $\mathcal{C}^{1+\varepsilon},$ has topological entropy, and the lift of $f$ to every finite covering of the torus is transitive, then the rotation set of $\widetilde f$ has non-empty interior.

Here we study homeomorphism of the 2-torus isotopic to the identity which have the following property:

\begin{definition}
A lift $\widetilde{f}$ of a homeomorphism of the torus isotopic to the identity is said to have the Property T if there exists a point $x_T$ in $\R{2}$ such that the omega-limit of the orbit of $x_T$ by $\widetilde{f},$ $\omega(x_T),$ contains all $\Z^2$ translates of $x_T,$ that is, $\pi^{-1}(\pi(x_T))\subset\omega(x_T).$ 
\end{definition}

The above hypothesis is clearly weaker than requesting transitivity of $\widetilde{f}.$ In fact, it covers the following case, which is of relevance for measure-preserving homeomorphisms, a hypothesis we do not require. A common and stable trait of sufficiently differentiable area preserving diffeomorphisms of the torus is the existence of elliptic islands. Of course such diffeomorphisms cannot be transitive, nor can their lift, but it is still possible for such mappings to be transitive in the complement of the interior of these islands, which in many relevant cases is a torus with infinitely many holes. If a diffeomorphism $f$ has a lift $\widetilde{f}$ which is transitive in the complement of the lift of the elliptic islands, then $\widetilde{f}$ obeys Property T.

Our main result is:

\begin{theorem}\label{teoprincipal}
Let $f$ be a homeomorphism of the 2-torus isotopic to the identity such that its lift $\widetilde{f}$ obeys Property T. Then $(0,0)$ is in the interior of $\rho(\widetilde{f}).$  
\end{theorem}

In section 2 we introduce the sets $B_{\theta}$ and $\omega(B_\theta),$ and we study their structure. In section 3 we study the projection of the sets $B_{\theta}$ and $\omega(B_{\theta})$ to the torus, and derive some properties of its closure. The final section is devoted to proving the main theorem.

In what follows $f$ will always be a homeomorphism of the 2-torus isotopic to the identity, and $\widetilde{f}$ a lift that satisfies Property T. $x_T$ will denote a point in $\R{2}$ such that $\pi^{-1}(\pi(x_T))=x_T+\Z^2 \subset\omega(x_T).$

\section{The sets $B_{\theta}, C_{\theta},\omega(B_{\theta}),$ and $\omega(C_{\theta})$}

 We begin by considering the one point compactification of $\R{2},$ the set $\R{2}\cup\{\infty\}$ which is homeomorphic to $S^2.$ We associate to the homeomorphism $f$ of the torus and its lift $\widetilde{f}$ to the plane a third homeomorphism $\widehat{f},$ a homeomorphism of $S^2$ induced by $\widetilde{f}$ that keeps the infinity fixed.

Given any $\theta$ in $[0, 2\pi),$ we define the sets
$$V_{\theta}=\{x)\in \R{2}: \langle x; e_\theta \rangle = 0\}$$
and 
 $$V^{+}_{\theta}=\{x\in \R{2}: \langle x; e_\theta \rangle \ge 0\},$$
where $e_{\theta} = (cos(\theta), sin(\theta)).$

We consider also the corresponding sets in $S^2$, which we denote by $\widehat{V_{\theta}}$ and $\widehat{V^{+}_{\theta}}.$
 
Let $\widehat{B_{\theta}}$ be the connected component of
$$\bigcap_{n\leq 0} \widehat{f}^n(\widehat{V}^{+}_{\theta}),$$
that contains the infinity, and denote by $B_{\theta}$ the corresponding set in $\R{2}.$ 

Let $\widehat{C_{\theta}}$ be the connected component of
$${\bigcap_{n\geq o} }\widehat{f}^n(\widehat{V}^{+}_{\theta}),$$
that contains the infinity, and denote by $C_{\theta}$ the corresponding set in $\R{2}.$ 

It follows from these definitions that $B_{\theta}$ is the union of all closed unbounded connected sets whose forward orbit remains in $V^{+}_{\theta}$ for all times, and that $C_{\theta}$ is the union of the connected components whose backward orbit remains in $V^{+}_{\theta}.$

\begin{proposition}\label{propbteta}
For every $\theta$ in $[0, 2\pi)$, we have:

\begin{enumerate}
\item{$B_{\theta}$ and $C_{\theta}$ are closed and $\widetilde {f}(B_{\theta})\subset B_{\theta}, \widetilde {f}^{-1}(C_{\theta})\subset C_{\theta}.$}
\item{Every connected component of $B_{\theta}$ or $C_{\theta}$  is unbounded.}
\item{If $(p,q)$ is a point in $V^{+}_{\theta}\cap \Z{}^2,$ then $B_{\theta}+(p,q)\subset B_{\theta},$ and the same holds for $C_{\theta}.$}
\item{$B_{\theta}^C$ has a single connected component, as does $C_{\theta}^C.$}
\end{enumerate}
\end{proposition}

\proof
 We prove this properties only for $B_{\theta}$, the reasoning is the same for $C_{\theta}.$ Properties 1 and 2 are direct consequences of the definition of $B_{\theta}$. If $x$ is a point of $B_{\theta},$ let $\Gamma$ be the connected component of $B_{\theta}$ that contains $x.$  Then, for every $y$ in $\Gamma$ and every positive integer $j,$ $\widetilde{f}^j(y)$ is in $V^{+}_{\theta}$. But this implies that, if $z$ is a point of the set $\Gamma + (p,q),$ then also for every positive integer $j$ we have $\widetilde{f}^j(z)\in V^{+}_{\theta}.$ Therefore $\Gamma + (p,q)$ is a connected closed set whose future orbit does not leave $V^{+}_{\theta}$, and so $\Gamma +(p,q)\subset B_{\theta},$ and $x +(p,q) \in B_{\theta}.$ This proves property 3.

Now let $\Omega$ be the connected component of $B_{\theta}^C$ that contains $\R{2}\setminus V^{+}_{\theta}.$ Since $\widetilde{f}(B_{\theta})\subset B_{\theta},$ we have that $B_{\theta}^C\subset \widetilde{f}(B_{\theta}^C),$ and so every connected component of $ \widetilde{f}(B_{\theta}^C)$ contains a connected component of $B_{\theta}^C.$ As $\widetilde{f}(\Omega)\cap\Omega\not=\emptyset,$ this implies $\Omega\subset\widetilde{f}(\Omega).$ In particular, $\Omega^C$ is forward invariant.

Let $z\in \Omega^C,$ and let $\Gamma$ be the closed connected component of $\Omega^C$ that contains $z$. Since $\Omega^C$ is forward invariant, for every positive integer $j$ we have $\widetilde {f}^j(\Gamma)\cap\Omega=\emptyset$ but since $(V^{+}_{\theta})^C\subset \Omega,$ $\Gamma$ is a closed connected set whose forward orbit remains in $V^{+}_{\theta}.$ Furthermore, since $\partial (\Omega^C)\subset B_{\theta},$ then $\Gamma\cap B_{\theta}\not=\emptyset.$ Let $z\in \Gamma\cap B_{\theta}$ and $\Theta$ be the connected component of $B_{\theta}$ that contains $z$. $\Theta$ is unbounded and since $\Theta \subset \Omega^C$ and $\Theta\cap\Gamma\not=\emptyset,$ then $\Theta\subset\Gamma.$ So  $\Gamma$ is unbounded, closed and its forward orbit remains in $V^{+}_{\theta},$ therefore $\Gamma\subset B_{\theta}.$ This shows that $\Omega^C=B_{\theta}$ and implies property 4 \endproof

\begin{lemma}
For every $\theta$ in $[0, 2\pi)$, we have $B_{\theta}\cap[-1,1]\times[-1,1]\not=\emptyset,$ and $C_{\theta}\cap[-1,1]\times[-1,1]\not=\emptyset.$
\end{lemma}
\proof We assume, without loss of generality, that $0\le \theta <\frac  {\pi}{2}.$  

Since $\widetilde {f}$ satisfies Property T, there is $y$ a translate of $x_T$ which lies in $V^{+}_{\theta+\pi}$ such that, for every positive integer $m,$ there exists a positive integer $n_0(m)$ such that $\widetilde{f}^{n_0(m)}(y)\in  V^{+}_{\theta}+(m,0).$ Let $n=n(m)$ be the first time that $\widetilde {f}^{n}(V^{+}_{\theta+\pi})$ intersects $V^{+}_{\theta}+(m,0)$, that is, 
$\widetilde {f}^{n}(V^{+}_{\theta+\pi})\cap V^{+}_{\theta}+(m,0)\not= \emptyset$ and such that, for all $0\le j<n,$ $\widetilde {f}^{j}(V^{+}_{\theta+\pi})\cap V^{+}_{\theta}+(m,0)$ is empty.

As $\widetilde{f}^{n(m)}(V^{+}_{\theta+\pi})\cap V^{+}_{\theta+\pi}$ is not empty, this implies that $$\left(\widetilde{f}^{n(m)}(\partial V^{+}_{\theta+\pi})\cap V^{+}_{\theta}+(m,0)\right) =\left(\widetilde{f}^{n(m)}(V_{\theta})\cap V^{+}_{\theta}+(m,0)\right)\not=\emptyset.$$

Let $m_k$ be a sequence of integers such that $\lim_{k\to\infty}n(m_k)=\infty.$ For each $k$ there is an unbounded closed curve $\gamma_k$ that starts at $\widetilde {f}^{n(m_k)}(V_{\theta})$ and such that
$\gamma_k\subset \widetilde {f}^{n(m_k)}(V^{+}_{\theta})\cap V^{+}_{\theta}+(m,0).$ It follows that the sets $\beta_k=\widetilde {f}^{-n(m_k)}(\gamma_k)$ all intersect $V_{\theta}$, and that for $0\le i\le n(m_k), \widetilde {f}^i(\beta_k)\subset V^{+}_{\theta}.$ 

Also, since for every point $z$ in $V_{\theta}$ there are integers $p$ and $q$ such that $ \langle (p,q);e_{\theta} \rangle \ge 0$ and such that $z+(p,q)$ is in $[0,1]\times[0,1],$ we have that there is  a curve $\alpha_k=\beta_k+(p,q)$ that starts in $[0,1]\times[0,1]$ and such that, for $0\le i\le n(m_k), \widetilde {f}^i(\alpha_k)\subset V^{+}_{\theta}.$

Now consider the closed sets $\widehat{\alpha_k}$ in $S^2$ corresponding to the sets $\alpha_k.$ $\widehat{\alpha_k}$ is a sequence of compact connected sets, and as such there is a subsequence $\widehat{\alpha_{k_j}}$ that is convergent in the Hausdorff topology to a closed connected set $\widehat{\Gamma}.$ It is clear that $\widehat{\alpha_k}$ is in $\widehat{B}_{\theta},$ and that $\widehat{\Gamma}$ intersects $[0,1]\times[0,1].$ The subset $\Gamma$ of $\R{2}$ corresponding to $\widehat{\Gamma}$ shows the lemma is true \endproof

Another interesting property of the set $B_{\theta}$ is the following:

\begin{lemma}
There exists $(p,q)\in \Z^{2}$ such that $\widetilde{f}^{-1}(B_{\theta})\subset B_{\theta}-(p,q).$ In particular, $f(\pi(B_{\theta}))=\pi(B_{\theta}).$
\end{lemma}

\proof Since $\widetilde{f}$ is biperiodic, there exists a positive real $K$ such that, for all $x\in\R{2}$, $\nrm{\widetilde{f}(x) - x}\le K.$

Let $y\in B_{\theta},$ and let $(p,q)\in\Z^{2}$ be such that $\langle (p,q); e_{\theta} \rangle \ge 2K.$ Let $\Gamma$ be the connected component of $y.$ Then $\widetilde{f}^{-1}(\Gamma) + (p,q)\subset V^{+}_{\theta}.$ As $\widetilde{f}(\widetilde{f}^{-1}(\Gamma) + (p,q))$ is a subset of $B_{\theta}+(p,q)\subset B_{\theta}\subset V^{+}_{\theta}$ and as $\widetilde{f}^{-1}(\Gamma) + (p,q)$ is connected and unbounded, we have that $\widetilde{f}^{-1}(\Gamma) + (p,q)\subset B_{\theta},$ which shows the result \endproof

Also, 
\begin{lemma}
There exists $(p,q)\in \Z^{2}$ such that $\widetilde{f}(C_{\theta})\subset C_{\theta}-(p,q).$ In particular, $f(\pi(C_{\theta}))=\pi(C_{\theta}).$
\end{lemma}
 and the proof is the same of the previous lemma.

We will need to also consider the following sets:
$$\omega(B_{\theta})=\cap_{i=1}^{\infty}\overline{\cup_{j=i}^{\infty} \widetilde{f}^j(B_{\theta})}=\cap_{i=1}^{\infty}\widetilde{f}^i(B_{\theta}),$$
where the last equality comes from the positive invariance of $B_{\theta},$ and the respective set in $S^{2}$
$$\widehat{\omega(B_{\theta})}=\cap_{i=1}^{\infty} \widehat{f}^i(B_{\theta}).$$ We also define the sets
$$\omega(C_{\theta}) =\cap_{i=1}^{\infty}\widetilde{f}^{-i}(C_{\theta}),$$
and $\widehat{\omega(C_{\theta})}.$

The following properties of $\omega(B_{\theta})$ and $\widehat{\omega(B_{\theta})}$ follow directly from the definition and proposition \ref{propbteta}.

\begin{proposition}\label{propomega}
For every $\theta$ in $[0, 2\pi)$, we have:
\begin{enumerate}
\item{$\{\infty\}\subset\widehat{\omega(B_{\theta})},$ and $\widehat{\omega(B_{\theta})}$ is connected.}
\item{$\omega(B_{\theta})$ is closed and $\widetilde {f}(\omega(B_{\theta}))=\omega( B_{\theta}).$}
\item{Every connected component of $\omega(B_{\theta})$ is unbounded.}
\item{If $(p,q)$ is a point in $V^{+}_{\theta}\cap \Z{}^2,$ then $\omega(B_{\theta})+(p,q)\subset \omega(B_{\theta}).$}
\item{$\omega(B_{\theta})^C$ has a single connected component.}
\end{enumerate}
\end{proposition}

\begin{proposition}
$\omega(B_{\theta})=\omega(C_{\theta}).$
\end{proposition}

\proof
We show only that $\omega(B_{\theta})\subset\omega(C_{\theta}).$ Let $\Gamma$ be a connected component of $\omega(B_{\theta}).$ Since $\omega(B_{\theta})$ is invariant, we have $\widetilde{f}^{-k}(\Gamma)\subset\omega(B_{\theta})\subset V^{+}_{\theta},$ for all positive integers $k,$ which shows that $\Gamma\subset C_{\theta}.$ Since the choice of $\Gamma$ was arbitrary, we have $\omega(B_{\theta})\subset C_{\theta}.$  Again, we use the invariance of $\omega(B_{\theta})$ and get
$$\widetilde{f}^{i}(\omega(B_{\theta}))=\omega(B_{\theta})\subset C_{\theta},$$
so that $\omega(B_{\theta})\subset \widetilde{f}^{-i}(C_{\theta})$ for all integers $i,$ proving the proposition \endproof

The following lemma will be our main tool in establishing the existence of vectors in the rotation set of $\widetilde{f}.$

\begin{lemma}\label{se_omega_vazio}
If $\omega(B_{\theta})=\emptyset,$ then there exists a positive real number $a$ such that:
\begin{itemize}
\item{For all $x$ in $B_{\theta},$
$$\liminf_{n\to\infty} \frac{\langle \widetilde{f}^{n}(x);e_{\theta}\rangle}{n}\ge a.$$}
\item{For all $x$ in $C_{\theta},$
$$\liminf_{n\to\infty} \frac{\langle \widetilde{f}^{-n}(x);e_{\theta}\rangle}{n}\ge a.$$}
\end{itemize}
\end{lemma}

\proof We again assume, without loss of generality, that $\theta$ is in $[0, \pi/2),$ and we prove just the first item. We claim that there is a positive $n$ such that $\widetilde {f}^i(B_{\theta})$ is contained in $V^{+}_{\theta}+(1,0),$ for all $i>n.$ If this was not the case, then there would be a sequence $i_j\to\infty$ such that, for all $j,$ there is a point $x_j=B_{\theta}$ such that $\widetilde {f}^{i_j}(x_j)\in V^{+}_{\theta}\setminus V^{+}_{\theta}+(1,0)$.

But this implies that there is a sequence of pair of integers $p_j, q_j$ such that
$\langle (p_j,q_j);e_{\theta} \rangle\ge 0$ and such that $\widetilde {f}^{i_j}(x_j)+(p_j,q_j)$ is in $[1,2]\times[0,1].$ But since $B_{\theta}+(p_j,q_j)\subset B_{\theta},$ we would have that, if $z_j=x_j+(p_j,q_j),$ then the sequence $(\widetilde {f}(z_j))_{j\in\N}$ would have a convergent subsequence, which contradicts $\omega(B_{\theta})=\emptyset,$ since $z_j \in B_{\theta}$ for all $j.$

Let $\Gamma$ be a connected component of $B_{\theta}.$ By the choice of $n$, for all $j\ge 0$ we have $\widetilde{f}^{j+n}(\Gamma)\subset V^{+}_{\theta}+(1,0),$ and so $\widetilde{f}^{n+j}(\Gamma)-(1,0)\subset V^{+}_{\theta}.$  Since $\widetilde{f}^n(\Gamma)$ is closed, connected and unbounded, it follows from the definition of $B_{\theta}$ that $\widetilde{f}^n(\Gamma)-(1,0)\subset B_{\theta},$  and so $\widetilde{f}^n(\Gamma)\subset B_{\theta}+(1,0).$ As the choice of $\Gamma$ was arbitrary, we have
$$\widetilde{f}^n(B_{\theta})\subset B_{\theta}.$$


A simple induction argument now shows that, for every positive integer $k$, $\widetilde {f}^{nk}(B_{\theta})\subset B_{\theta}+(k,0)\subset V^{+}_{\theta}+(k,0),$ and so, for every $x$ in $B_{\theta}$,
$$
\liminf_{i\to\infty}\frac{\langle \widetilde {f}^i(x);e_{\theta}\rangle}{i}\ge \frac{cos(\theta)}{n},
$$
which ends the proof \endproof 

The following is an immediate consequence of the previous lemma:
\begin{corollary}\label{vetorrotacaonosdoisplanos}
If $\omega(B_{\theta})=\emptyset,$ than there exists $r_1, r_2\in\rho(\widetilde{f})$ such that $\langle r_1;e_\theta\rangle>0,$ and $\langle r_2;e_\theta\rangle<0.$
\end{corollary}  

We can now state the main lemma of this work: which will be proved in the last section:

\begin{lemma}\label{mainlemma}
For every $\theta\in[0,\pi),$ at least one of $\omega(B_{\theta})$ and $\omega(B_{\theta+\pi})$ is empty.
\end{lemma}

We will postpone the proof of this result to the last section, but let us now show how this lemma yields theorem \ref{teoprincipal}.

\

{\noindent \bf Proof of theorem \ref{teoprincipal}:}
We first note that, since $\widetilde{f}$ satisfies Property T, it has a recurrent point, and as so $(0,0)\in \rho(\widetilde{f}).$ Since $\rho(\widetilde{f})$ is convex, either $(0,0)$ is in the interior of $\rho(\widetilde{f}),$ and we are done, or there exists a straight line passing through $(0,0)$ such that $\rho(\widetilde{f})$ is contained in one side of this line. Assume by contradiction that the latter happens. This means that there is $\theta_1$ such that for all $r\in \rho(\widetilde{f}),$ $\langle r;e_{\theta_1}\rangle \geq 0.$

But from lemma \ref{mainlemma} either $\omega(B_{\theta_1}) = \emptyset$ or $\omega(B_{\theta_1 + \pi})=\emptyset,$ and in both cases this implies, by corollary \ref{vetorrotacaonosdoisplanos}, that there exists $r\in\rho(\widetilde{f})$ with $\langle r;e_{\theta_1}\rangle < 0.$ \endproof

\section{The sets $\pi(B_{\theta})$ and $\pi(\omega(B_{\theta}))$}

Let us state this general result:

\begin{lemma}
Let $g$ be a homeomorphism of the torus isotopic to the identity and $\widetilde g$ be its lift to the plane. Let $\widetilde A\subset\R{2}$ be such that if $\widetilde  z \in \widetilde A$, then $\widetilde  z + (i,0) \in \widetilde A$ for all integer $i.$ Assume that there exists a point $\widetilde  x$ in $\R{2}$ such that 
$$\widetilde{A}\subset\omega(\widetilde{x})=\cap_{i=1}^{\infty}\overline{\cup_{j=i}^{\infty} \widetilde{g}^j(\widetilde{x})}.$$

Then, for every $z\in A=\pi(\widetilde A)$ and every $\varepsilon>0$, there exists a curve $\gamma\subset \R{2}$ which is unbounded in the $e_0$ direction but bounded in the $e_{\theta/2}$ direction, and such that, for every point in $y\in\gamma$ there exists a positive $n(y)$ such that $\pi(\widetilde {g}^{n(y)}(y))\in B_{\varepsilon}(z)$.
\end{lemma}

{\it Proof:}

Let $z\in A$ and and $\widetilde  z$ be a point such that $\pi(\widetilde  z)=z.$ Since $\widetilde{A}\subset\omega(\widetilde{x}),$ there are points $x_1$ and $x_2$ and integers $k_1$ and $k_2$, both strictly negative, such that 
$\widetilde  g ^{k_1}(x_1) = y_1 \in B_{\varepsilon}(\widetilde  z)$ and $\widetilde g ^{k_2}(x_2) = y_2 \in B_{\varepsilon}(\widetilde  z + (1,0)).$ Let $\alpha_0$ be the segment connecting $x_1$ and $x_2,$ $\alpha_1$ be the segment connecting $x_1$ and $y_1$, $\alpha_2$ be the segment connecting $x_2$ and $y_2 - (1,0),$   $\alpha_3$ be the segment connecting $y_1$ and $\widetilde  z$ and $\alpha_4$ be the segment connecting $y_2-(1,0)$ and $\widetilde  z$. We note that $\pi(\alpha_i)\subset B_{\varepsilon}(z),$ for $i\in \{0,1,2,3,4\}.$

Also, since $x_1$ and $y_1$ belong to $\alpha_1,$ and since $y_1=\widetilde{g}^{k_1}(x_1)$ belongs to $\widetilde{g}^{k_1}(\alpha_1),$ we have that $\widetilde{g}^{k_1}(\alpha_1)\cup\alpha_1$ is a connected set that contains $\widetilde{g}^{k_1}(x_1).$ A simple induction shows that 
$$\beta_1=\left(\bigcup_{i=0}^{-k_2}\widetilde{g}^{ik_1}(\alpha_1)\right),$$
is a connected arc joining $\widetilde{g}^{-k_1k_2}(x_1)$ and $y_1$.

Analogously, since $x_2$ and $y_2 - (1,0)$ belong to $\alpha_2$ and $\widetilde{g}^{k_2}(x_2)=y_2,$ it follows that the set $\widetilde{g}^{k_2}(\alpha_2)\cup\alpha_2+(1,0)$ is an arc connecting $\widetilde{g}^{k_2}(x_2)$ and $y_2.$ Therefore the set
$$\beta_2=  \left(\bigcup_{i=0}^{-k_1}(\widetilde{g}^{ik_2}(\alpha_2+(i-k_1,0))\right)
$$
is a connected arc joining $\widetilde{g}^{-k_1k_2}(x_2)$ and $y_2-(k_1-1,0)$.
This implies that 
$$\widetilde {g}^{-k_1k_2}(\alpha_0)\cup\beta_1
\cup\beta_2$$
is a connected arc joining $y_1$ and $y_2-(k_1-1,0).$ 

Consider the set
$$\gamma_0=\widetilde {g}^{-k_1k_2}(\alpha_0)\cup\beta_1
\cup\beta_2\cup\alpha_3\cup\alpha_4+(-k_1,0),$$
which is also connected and bounded in the $e_{\pi/2}$ direction, and which contains both $\widetilde  z$ and $\widetilde  z -(k_1,0).$ If $\gamma=\bigcup_{i=-\infty}^{\infty}(\gamma_0 + (ik_1,0))$, then $\gamma $ is connected, unbounded in the $e_0$ direction and bounded in the $e_{\pi/2}$ direction, and satisfies that, if $y\in\gamma$, then there is a $n(y)>0$ such that $\pi(\widetilde{g}^{n(y)}(y))\in \pi(\alpha_i)$ for some $0\le i\le 4.$ Since $\pi(\alpha_i)\subset B_{\varepsilon}(z), 0\le i\le 4,$ we have the result \endproof

The previous lemma yields this two simple corollaries, the first one which will be used in the remainder of the paper, the second interesting in itself.

\begin{corollary}\label{denso em tildeA}
Let $g$ be a homeomorphism of the torus and $\widetilde  g$ be its lift to the plane. Assume there is a set $\widetilde  A\subset \R{2}$ which is $\Z{}^2$ invariant, i.e., for every $\widetilde {y}$ in $\widetilde  A,$ and every pair of integers $(p,q)\in \Z{2},$ $\widetilde  y +(p,q)\in \widetilde  A.$ Assume further than there is a point $\widetilde  x$ in $\widetilde  A$ such that $\widetilde A\subset \omega(\widetilde x).$ 

Let $B\subset \R{2}$ be a set such that $\widetilde {g}(B)\subset B$ and such that $B$ has an unbounded connected component. Then $\pi(\widetilde  A)\subset \overline{\pi(B)}.$
\end{corollary}

\proof Let $\widetilde {z}\in B$ be a point in an unbounded connected component of $B.$
Suppose, by contradiction, that there exists a point $\widetilde {y}\in\widetilde {A}$ and a positive real $\varepsilon$ such that $\pi(B)\cap B_{\varepsilon}(\pi(\widetilde  y))=\emptyset.$

The previous lemma implies that there are two connected sets, $\alpha$ and $\beta,$ with the following property:
\begin{itemize}
\item{For all $\widetilde {x}$ in $\alpha\cup\beta$, there exists a positive integer $n(\widetilde x)$ such that $\pi(\widetilde {g}^n(\widetilde {x}))\in  B_{\varepsilon}(\pi(\widetilde  y)).$}
\item{ $\alpha$ is unbounded in the direction of $e_{\pi/2},$ bounded in the direction of $e_{0}$ and disconnects the plane.}
\item{ $\beta$ is bounded in the direction of $e_{\pi/2},$ unbounded in the direction of $e_{0}$ and disconnects the plane.}
\end{itemize}

Now, by the contradiction hypothesis and the first item, 
\begin{eqnarray}\label{aux1}
\pi(B)\cap\pi(\alpha\cup\beta)=\emptyset
\end{eqnarray}
By the second item, there is an integer $k_1$ such that $\widetilde {z}$ is to the left of $\alpha+( k_1, 0),$  to the right of $\alpha-(k_1, 0),$  below $\beta + (0, k_1)$ and above $\beta -(0, k_1).$ This shows that $\widetilde {z}$ is in a bounded connected component of 
the complement 
 $$D=\alpha+( k_1, 0)\cup\alpha-(k_1, 0)\cup\beta + (0, k_1)\cup  \beta -(0, k_1).$$ 

But, from (\ref{aux1}),  $B\cap D=\emptyset,$ which is absurd since the connected component of $B$ that contains $\widetilde {z}$ must be unbounded by hypothesis and remain in the same connected component of $D^C$ that contains $\widetilde {z}.$ \endproof

We finish this section with 2 direct applications of this results:
\begin{corollary}
Let $g$ be a homeomorphism of the torus isotopic to the identity and $\widetilde  g$ be its lift. Assume $\widetilde  g$ is transitive in $\R{2}.$  

Let $B\subset \R{2}$ be a set such that $\widetilde {g}(B)\subset B$ and such that $B$ has an unbounded connected component. Then $\overline{\pi(B)}=T^2.$
\end{corollary}

\begin{corollary}\label{xT_tanofecho}
Let $\widetilde f$ be a homeomorphism of the torus satisfying property T. Then, for all $0\le\theta<2\pi,$ if $\omega(B_{\theta})$ is non-empty, then
$\pi(x_T)\subset \overline \pi(\omega(B_{\theta})).$
\end{corollary}

\section{Proof of lemma \ref{mainlemma}}

It suffices to show that, for every $\theta\in[0,\pi),$ either $\omega(B_{\theta})=\emptyset$, or $\omega(B_{\theta+\pi})=\emptyset.$

Assume, by contradiction, that there exists $\theta$ for which this does not hold. Without loss in generality, we will assume that $0\le \theta<\pi/2.$

There are 2 possibilities:
\begin{enumerate}
\item{$\pi(\omega(B_\theta))\cap\pi(\omega(B_{\theta+\pi}))\not=\emptyset$}
\item{$\pi(\omega(B_\theta))\cap\pi(\omega(B_{\theta+\pi}))=\emptyset$}
\end{enumerate}

\subsection{ Case 1}

If $\pi(\omega(B_\theta))\cap\pi(\omega(B_{\theta+\pi}))\not=\emptyset,$ then there exists $(p_1,q_1)$ and $(p_2,q_2)$ such that
$$\omega(B_{\theta})+(p_1,q_1)\cap \omega(B_{\theta+\pi})+(p_2,q_2)\not=\emptyset,$$
and so, if $(p,q)=(p_2-p_1,q_2-q_1),$
$$\omega(B_{\theta})\cap \omega(B_{\theta+\pi})+(p,q)\not=\emptyset.$$

In this case, we have:
\begin{lemma}\label{omegaintersecta}
The open set 
$$ O_1=\left((\omega(B_{\theta})-(1,1))\bigcup (\omega(B_{\theta+\pi})+(p+1,q+1))\right)^C,$$
has infinitely many connected components.
Furthermore, the set $x_T+\Z^2$ intersects infinitely many of these connected components.
\end{lemma}

\proof
Let $N$ be a fixed integer. There exists integers
$n_1,m_1, n_2, m_2$ such that $n_2<0<n_1, m_1>0, m_2>0$ and such that
$${\frac{-1}{4N}}< \langle(n_1,m_1);e_{\theta}\rangle <0,$$
and 
$$0<\langle (n_2,m_2);e_{\theta}\rangle < {\frac{1}{4N}}.$$

Let $y\in x_T+\Z^2$ be such that
$1/4\le \langle y; e_{\theta}\rangle\le 3/4.$ 

We will show that at least $N$ of the points $y, y+i(n_1,m_1), y+j(n_2,m_2)$, for $1\le i\le N, 1\le j\le N,$ are in distinct connected components of 
$ O_1,$
which clearly proves the proposition.

Assume, by contradiction, that there are $i_1,i_2, j_1,$ and $j_2$ in $[0, N],$ $i_1<i_2,$ $j_1<j_2$ such that $y_1=y+i_1(n_1,m_1)$ and $y_1 + (i_2-i_1)(n_1,m_1)$ are in the same connected component, as are $y_2=y+j_1(n_2,m_2)$ and $y_2+(j_2-j_1)(n_2,m_2).$ 

Since $y_1$ and $y_1 + (i_2-i_1)(n_1,m_1)$ are in the same connected component, there is a curve $\alpha$ joining these 2 points that does not intersect $O_1^C.$
Likewise, there is a curve $\beta$ joining $y_2$ and $y_2+(j_2-j_1)(n_2,m_2),$ that also does not intersect $O_1^C.$ Note that 
\begin{eqnarray}\label{0<y<1}
0<\langle y_1;e_{\theta}\rangle<1,\\
0<\langle y_2;e_{\theta}\rangle<1. \nonumber
\end{eqnarray}

We recall that, if $(r_1,s_1)\in\Z^2$ is such that 
$\langle (r_1,s_1);e_{\theta}\rangle < 1 \le \langle (1,1);e_{\theta}\rangle,$
then $\langle (1-r_1,1-s_1);e_{\theta}\rangle >0,$ and from proposition \ref{propomega} we have $\omega(B_{\theta})+(1-r_1,1-s_1)\subset \omega(B_{\theta}),$ and $\omega(B_{\theta})\subset \omega(B_{\theta})+(r_1-1,s_1-1).$ Therefore, if $C\cap\omega(B_{\theta})-(1,1)=\emptyset,$ than 
$C+(r_1,s_1)\cap\omega(B_{\theta})+(r_1-1,s_1-1)=\emptyset$ and so $C+(r_1,s_1)\cap\omega(B_{\theta})=\emptyset.$

Likewise, if $(r_2,s_2)\in \Z^2$ is such that $\langle (r_2,s_2);e_{\theta}\rangle > -1,$ and if $C\cap \omega(B_{\theta+\pi})+(p+1,q+1)=\emptyset,$ then $C+(r_2,s_2)\cap \omega(B_{\theta +\pi}) +(p,q)=\emptyset.$

This and (\ref{0<y<1}) shows that, if $z$ is such that $\pi(z)=\pi(y)$ and such that $0\le \langle z-y;e_{\theta}\rangle\le 1,$ then both $\alpha+(z-y)$ and $\beta+(z-y)$ do not intersect $\omega(B_{\theta})\cap \omega(B_{\theta+\pi})+(p,q).$ 

Let $z_0\in\pi^{-1}(\pi(x_T))$ be such that $1/4<\langle z_0;e_\theta\rangle<3/4.$ For $-\infty\le i\le +\infty,$ we define, if $\langle z_i;e_\theta\rangle<1/2,$
$\gamma_i= \alpha+(z_i-y_1), z_{i+1}= z_i + (i_2-i_1)(n_1,m_1);$  and if $\langle z_i;e_\theta\rangle>1/2,$ we define $\gamma_i= \beta+(z_i-y_2), z_{i+1}= z_i + (j_2-j_1)(n_2,m_2).$ Note that $\gamma_i$ connects $z_i$ and $z_{i+1}.$  Also notice that, since $m_1$ and $m_2$ are positive, $\langle z_{i+1}, e_{\pi/2}\rangle > \langle z_i, e_{\pi/2}\rangle$ and that
$$\lim_{i\to\infty} \langle z_i, e_{\pi/2} \rangle = +\infty, \lim_{i\to-\infty} \langle z_i, e_{\pi/2} \rangle = -\infty.$$

Also, since $-\frac{1}{4}\langle (i_2-i_1)(n_1,m_1);e_{\theta}\rangle<0$ and $0<\langle (j_2-j_1)(n_2,m_2);e_{\theta}\rangle<\frac{1}{4},$ then an induction shows that $1/4<\langle z_i;e_\theta\rangle<3/4,$ for all integers $i.$ Therefore, for all $i\in\Z,$ both $\abs{\langle z_i-y_1;e_{\theta}\rangle}<1$ and 
$\abs{\langle z_i-y_2;e_{\theta}\rangle}<1.$ This shows that $\gamma_i\cap\left(\omega(B_{\theta})\cup \omega(B_{\theta+\pi})+(p,q)\right)=\emptyset.$

Finally, let $M=\max_{x\in(\alpha\cup\beta)}\{\abs{\langle x;e_\theta\rangle}\}.$ Then 
$$\max_{x\in\gamma_i}\{\abs{\langle x;e_\theta\rangle}\}<M+1.$$ 

Let $\gamma=\cup_{1=-\infty}^{\infty}\gamma_i.$ Then clearly $\gamma$ is connected, does not intersect $\omega(B_{\theta})\cup \omega(B_{\theta+\pi})+(p,q),$ is unbounded in the $e_{\pi/2}$ direction, and bounded in the $e_{\theta}$ direction. This shows that the sets $\{x\in\R{2}\mid \langle x, e_{\theta}\rangle > M+1\},$ and $\{x\in\R{2}\mid \langle x, e_{\theta}\rangle < -(M+1)\}$ are in different connected components of $\gamma^C.$

But since $\omega(B_{\theta})\cap \omega(B_{\theta+\pi})+(p,q)$ is not empty, there exists $x$ in the intersection. Let $\Gamma$ be the connected component of $\omega(B_{\theta})$ that contains $x,$ and let $\Theta$ be the connected component of $\omega(B_{\theta+\pi})+(p,q)$  that contains $x.$ Then $\Gamma_1= \Gamma \cup \Theta$ is a connected component of $\omega(B_{\theta})\cap \omega(B_{\theta+\pi})+(p,q)$ such that $sup_{w\in\Gamma_1}\langle w;e_{\theta})\rangle=+\infty$ and such that
$inf_{w\in\Gamma_1}\langle w;e_{\theta}\rangle=-\infty.$ But this contradicts the fact that $\gamma\cap\Gamma_1=\emptyset$ and that $\gamma$ is connected \endproof

Since $\widetilde {f}$ satisfies Property T, the forward orbit of $x_T$ is not confined to any half-plane and, as such, 
$x_T\notin(\omega(B_{\theta})\cup\omega(B_{\theta+\pi})+(p,q),$ so that $x_T\in O_1.$

Now $x_T$ is a recurrent point, so there is an $n$ such that $\widetilde {f}^n(x_T)$ and $x_T$ are in the same connected component of
$O_1.$ Therefore, if $U$ is the connected component of $O_1$ that contains $x_T,$ then $\widetilde {f}^n(U)\cap U$ is not empty, and since $O_1$ is an invariant set, this implies that $\widetilde {f}^n(U)=U.$ Thus, the orbit of $x_T$ only intersects $n$ connected components of $O_1.$ 

But lemma \ref{omegaintersecta} tells us that the set $x_T+\Z^2$ both intersects infinitely many connected components of  $O_1,$ and is contained in the closure of the forward orbit of $x_T$ which is absurd. Therefore Case 1 cannot be.

\subsection{ Case 2}

We still need to deal with the case where
$$\pi(\omega(B_{\theta}))\cap \pi(\omega(B_{\theta +\pi}))=\emptyset,$$
which we assume throughout the section.

We begin with the following proposition, whose proof can be seen in \cite{plane}

\begin{proposition}
Let $K_1$ and $K_2$ be compact connected sets of $\R{2}$ such that $K_1\cap K_2$ is connected. If $x$ and $y$ are in the same connected component of $K_1^C$ and are also in the same connected component of $K_2^C,$ then $x$ and $y$ are in the same connected component of $(K_1\cup K_2)^C.$
\end{proposition}

\begin{proposition}\label{temso1comp}
For every $(p_1,q_1), (p_2,q_2)\in \Z^2,$ the complement of the set 
$$\left( \omega(B_{\theta})+(p_1,q_1)\cup \omega(B_{\theta+\pi})+(p_2,q_2)\right)$$
has a single connected component.
\end{proposition}

\proof We already know that $\omega(B_{\theta})+(p_1,q_1)\cap \omega(B_{\theta+\pi})+(p_2,q_2)=\emptyset.$ We also know that, from lemma \ref{propbteta}, that both the complements of $\omega(B_{\theta})+(p_1,q_1)$ and of $\omega(B_{\theta+\pi})+(p_2,q_2)$ have a single connected component. This means that, if we consider the corresponding sets in $S^2,$
$\widehat{\omega(B_{\theta})+(p_1,q_1)}$ and $\widehat{ \omega(B_{\theta+\pi})+(p_2,q_2)},$ their intersection consists of the point $\infty,$ and so we can apply the previous proposition, which gives us the result \endproof

The following is result follows from a construction very similar to the one done in the proof of lemma \ref{omegaintersecta}.

\begin{proposition}\label{separaomegas}
For every point $x$ in $T^2,$ and every $(p_1,q_1)$ and $(p_2, q_2)$ in $\Z^2$, there is an arc $\gamma\subset\R{2}$ that contains infinitely many points of $\pi^{-1}(x)$ and such that $\omega(B_{\theta})+(p_1,q_1)$ and $\omega(B_{\theta+\pi})+(p_2,q_2)$ are in distinct connected components of $\gamma^C.$
\end{proposition} 

\proof Since $\pi(\omega(B_{\theta}))\cap \pi(\omega(B_{\theta} +\pi))=\emptyset,$ we will assume, without loss in generality, that $x\notin \pi(\omega(B_{\theta})).$ Therefore, there is a point $y$ in $\pi^{-1}(x)$ such that $\langle y;e_\theta\rangle \ge \langle(p_2,q_2); e_\theta\rangle +3.$ This means that the points
$y$, $y+(0,1)$ and $y-(1,0)$ are all in the semi-plane $V^{+}_{\theta}+(p_2+1,q_2+1),$ and therefore none of these 3 points belongs to 
$\omega(B_{\theta})+(p_1-1,q_1-1) \cup \omega(B_{\theta+\pi})+(p_2+1,q_2+1).$

Now, from proposition \ref{temso1comp}, we know that all three points are in the same connected component of $\left(\omega(B_{\theta})+(p_1-1,q_1-1) \cup \omega(B_{\theta+\pi})+(p_2+1,q_2+1)\right)^C,$ and so there is an arc $\alpha$ joining $y$ and $y+(0,1)$ and an arc $\beta$ joining $y$ and $y-(1,0)$ such that
 $$(\alpha\cup\beta)\cap \left(\omega(B_{\theta})+(p_1-1,q_1-1) \cup \omega(B_{\theta+\pi})+(p_2+1,q_2+1)\right)=\emptyset.$$

Now an argument similar to the one as in the proof of lemma \ref{se_omega_vazio}, can construct the arc $\gamma$ by joining together translations of $\alpha$ and $\beta,$ such that 
$$-\infty<m= \inf_{z\in\gamma} \langle z; e_\theta\rangle \le \sup_{z\in\gamma} \langle z; e_\theta \rangle =M< \infty,$$
and such that $\gamma^C$ has two connected unbounded components, one of which contains all points of $V^{+}_{\theta}+M e_\theta,$ and the other contains all points of $\R{2}\setminus (V^{+}_{\theta}+m e_\theta).$

Since $\omega(B_{\theta})+(p_1,q_1)\cap V^{+}_{\theta}+M e_\theta\not=\emptyset,$ and since
$$\omega(B_{\theta+\pi})+(p_2,q_2)\cap \R{2}\setminus (V^{+}_{\theta}+m e_\theta)$$
 is also nonempty, we have the result\endproof  

Since $x_T+\Z^2\subset \omega(x_T),$ $\pi(x_T)$ is not a fixed point. Let $\varepsilon>0$ such that, if $B_{\varepsilon}[\pi(x_T)]$ is the closed ball with center $\pi(x_T)$ and radius $\varepsilon$, then $f(B_{\varepsilon}[\pi(x_T)])\cap B_{\varepsilon}[\pi(x_T)]=\emptyset.$

Now, from corollary \ref{xT_tanofecho}, $$\pi(x_T)\in\overline{\pi(\omega(B_{\theta}))}\cap \overline{\pi(\omega(B_{\theta} +\pi))},$$
therefore there exists $y\in x_T+\Z^2,$ and two integers $p$ and $q,$ such that  the distance from $y$ to $B_{\theta}$ is less than $\varepsilon$, as is the distance between $y-(p,q)$ to $B_{\theta+\pi}.$

Since we assumed that $\pi(\omega(B_\theta))\cap \pi(\omega(B_{\theta+\pi}))=\emptyset,$ the distance between the compact sets $\omega(B_\theta)\cap B_\varepsilon[y]$ and $(\omega(B_{\theta+\pi})+(p,q))\cap B_\varepsilon[y]$ is strictly positive number $a,$ and $a<\varepsilon.$ Let $v$ be an open line segment with length $a$ such that its endpoints are $u$ and $w,$ where $u\in \omega(B_\theta)$ and $w\in \omega(B_{\theta+\pi})+(p,q).$ Note that, as $v\subset B_\varepsilon[y], \widetilde{f}(V)\cap V= \emptyset.$

Finally, let $\Gamma$ be the connected component of $\omega(B_\theta)$ that contains $u$ and let $\Theta$ be the connected component of $\omega(B_{\theta+\pi})+(p,q)$ that contains $w$.

\begin{lemma}\label{tem2componentes}
The set $(\Gamma\cup\Theta)^C$ contains a single connected component. The set $(\Gamma\cup\Theta\cup v)^C$ has exactly two connected components, $\Omega_1$ and $\Omega_2$. Also, for every $x$ in $T^2,$ $\pi^{-1}(x)$ intersects both $\Omega_1$ and $\Omega_2.$

\end{lemma}

\proof The first claim follows from proposition \ref{temso1comp}.

Let $p$ be the midpoint of the segment $v,$ and let $r= dist (p, \Theta)= dist (p, \Gamma).$ The closed ball $B_r[p]$ is splitted in two by $v.$  Of course, since 
$(\Gamma\cup \Theta)^C$ has a single connected component,  $(\Gamma\cup\Theta\cup v)^C$ has, at most, two components. Let $p$ be a point in one of the connected components of  $B_r[p]\setminus v$ and let $q$ be a point in the second connected component. Suppose, by contradiction, that $p$ and $q$ are in the same connected component of $(\Gamma\cup\Theta\cup v)^C.$
Then, there is an arc $\alpha$ connecting $p$ and $q$ lying entirely on this set. We can add to $\alpha$ the segment joining $p$ and $q$ to obtain a Jordan curve $J$, which does not intersect $\Gamma\cup\Theta.$ But one endpoint of $v$ must be in the $int(J)$, and so one of $\Gamma$ and $\Theta$ must be in the $int(J),$ absurd since both are unbounded.  $\endproof$
 
\begin{proposition}\label{autointersect}
For $i=\{1,2\},\widetilde {f}(\Omega_i)\cap\Omega_i$ is not empty.
\end{proposition}
\proof Let $y$ be a fixed point of $\widetilde{f}$ (Since $x_T$ is a recurrent point of $\widetilde{f},$ $\widetilde{f}$ cannot be a Brouwer homeomorphism and so it has fixed points). By proposition \ref{separaomegas}, there is an arc $\gamma$ which connects infinitely many translates of $y$ and that separates  $\omega(B_\theta)$ and $\omega(B_{\theta+\pi})+(p,q).$ Since $\Gamma\cup\Theta\cup v$ is connected, $\gamma\cap v\not=\emptyset,$ and so there must be translates of $y$ both in $\Omega_1$ and in $\Omega_2$ \endproof 

\begin{lemma}
Either $\widetilde {f}(\overline{\Omega_1})\subset \overline{\Omega_1}$ or $\widetilde {f}(\overline{\Omega_2})\subset \overline{\Omega_2}.$ 
\end{lemma} 

\proof
By the definition of $v$, $v\cap\omega(B_\theta)=\emptyset,$ as well as $v\cap \omega(B_{\theta+\pi})+(p,q)=\emptyset.$ Since both $\omega(B_\theta)$ and $\omega(B_{\theta+\pi})$ are invariant, and since $\widetilde {f}(B_\varepsilon[y])\cap B_\varepsilon[y]=\emptyset,$ we have that
$\widetilde {f}(v)\cap(\Gamma\cup\Theta\cup v)=\emptyset.$

This means that either $\widetilde {f}(v)\subset \Omega_1$ or $\widetilde {f}(v)\subset \Omega_2.$ We will assume, without loss of generality, that the former is true.

Since $\Gamma$ is a connected component of $\omega(B_\theta)$ and the latter is invariant, then $\widetilde{f}(\Gamma)$ is also a connected component of $\omega(B_{\theta})$ and so either $\widetilde {f}(\Gamma)=\Gamma$ or $\widetilde {f}(\Gamma) \cap \Gamma=\emptyset.$
 We also have that $\widetilde {f}(\Gamma)\cap(v\cup\Theta)=\emptyset.$ If $\widetilde {f}(\Gamma)=\Gamma,$ then $\widetilde {f}(\Gamma)\cap \Omega_2=\emptyset.$ If $\widetilde {f}(\Gamma)\not=\Gamma$, since $\Gamma\cup v$ is connected, and since $\widetilde {f}(v)\subset \Omega_1$ and $\widetilde {f}(\Gamma)\cap\partial \Omega_1=\emptyset,$ we have $\widetilde {f}(\Gamma)\subset \Omega_1,$ and thus $\widetilde {f}(\Gamma)\cap \Omega_2$ is also empty.
 
By a similar argument, $\widetilde {f}(\Theta)\cap \Omega_2=\emptyset,$ and so $\widetilde{f}(\partial\Omega_1)\cap \Omega_2$ is empty. Since both $\Omega_1$ and $\Omega_2$ are connected, this implies that either $\Omega_2\subset\widetilde {f}(\Omega_1)$ or $\Omega_2\cap \widetilde{f}(\Omega_1)=\emptyset.$ By proposition \ref{autointersect} the first possibility cannot be, since $\widetilde{f}$ is a homeomorphism. The second possibility implies that and $\widetilde {f}(\overline{\Omega_1})\subset \overline{\Omega_1}.$ \endproof 

But since $\widetilde{f}$ satisfies property T the future orbit of $x_T$ must visit any neighborhood of any $\Z^2$ translate of $x_T$ infinitely many times, and so, from lemma \ref{tem2componentes}, there must be integers $n_1<n_2<n_3$ such that 

\begin{eqnarray}\label{n1n2n3}
\widetilde{f}^{n_1}(x_T)\in \Omega_1, \widetilde{f}^{n_2}(x_T)\in \Omega_2 \hbox{ and}\, \widetilde{f}^{n_3}(x_T)\in \Omega_1.
\end{eqnarray}

 But, from the previous lemma, if $\widetilde{f}(\overline{\Omega_1})\subset\overline{\Omega_1},$ then for all $n\ge n_1,\, \widetilde{f}^{n}(x_T)\notin \Omega_2$, which violates (\ref{n1n2n3}) and if $\widetilde{f}(\overline{\Omega_2})\subset\overline{\Omega_2},$  then for all $n\ge n_2,\, \widetilde{f}^{n}(x_T)\notin \Omega_1$, which again violates (\ref{n1n2n3}), and we are done.

\section{aknowledgment}

 I'd like to thank Prof. Addas-Zanata for the insightful discussions on this work, and presenting ideas that helped to clarify some of the proofs presented here.
 
\bibliographystyle{amsplain}

\end{document}